\documentclass[12pt,leqno]{article}
\usepackage{amsfonts}
\usepackage{amsmath, amsthm, amsfonts, amssymb, color,latexsym, enumerate}
\usepackage{mathrsfs}
\usepackage{url}
\usepackage{color}

\newtheorem{thm}{Theorem}[section]
\newtheorem{cor}[thm]{Corollary}
\newtheorem{lem}[thm]{Lemma}
\newtheorem{prp}[thm]{Proposition}

\newtheorem{rem}[thm]{Remark}
\newtheorem{asm}[thm]{Assumption}
\theoremstyle{definition}

\newcommand{\scr}[1]{\mathscr #1}
\definecolor{wco}{rgb}{0.5,0.2,0.3}

\numberwithin{equation}{section} \theoremstyle{remark}

\title{{\bf   Bernstein-type inequalities for Markov processes via resolvent estimates   }\footnote{Supported by the Natural Science Foundation of Shanghai (25ZR1402414) and the National Natural Science Foundation of China (12501185).} }
\author{
{\bf  Jie-Xiang Zhu$^{a)}$  }\\
\footnotesize{$^{a)}$ Department of Mathematics, Shanghai Normal University, Shanghai 200234, China }\\
 \footnotesize{  jiexiangzhu7@gmail.com } }

\begin{document}
\allowdisplaybreaks
\def\R{\mathbb R}
\def\C{\mathbb C}
\def\N{\mathbb N}
\def\<{\langle} \def\>{\rangle}
 \def\E{\mathbb E}
\def\d{\text{\rm{d}}}
\def\beg{\begin} \def\beq{\begin{equation}}
\def\e{\text{\rm{e}}}
 \def\tt{\tilde}
\def\cut{\text{\rm{cut}}} \def\P{\mathbb P} \def\ifn{I_n(f^{\bigotimes n})}
\def\aaa{\mathbf{r}}
\def\gap{\text{\rm{gap}}}
\def\Z{\mathbb Z}
\def\L{\mathcal L}
\def\M{\scr M} \def\Q{\mathbb Q}
\def\Rank{{\rm Rank}} \def\B{\scr B} \def\i{{\rm i}} \def\HR{\hat{\R}^d}
\def\to{\rightarrow}\def\l{\ell}\def\iint{\int}
\def\EE{\scr E}
\def\A{\scr A} \def\Lip{{\rm Lip}}
 \def\i{{\rm i}}
\def\PP{{\bf P}}\def\HL{{\bf L}}
\def\Ss{\mathbb S}
\def\BL{\scr A}\def\Pp{\mathbb P} \def\Ee{\mathcal E}
\def\V{\mathrm V}
\def\D{\mathscr D}

\maketitle

\begin{abstract}
We establish Bernstein-type inequalities for additive functionals of Markov processes through uniform resolvent estimates between a Banach space and its dual.
Under an iterated Poincar\'e inequality, we obtain the classical Bernstein form sought by Huang and Li in \cite{HL}, up to numerical constants, in both continuous and discrete time. Suitable choices of the Banach space recover Bernstein-type bounds under Poincar\'e and Sobolev inequalities and allow
unbounded observables. For normal generators, a criterion involving the modulus of the generator captures the contribution of its antisymmetric part.
As an application, we obtain Bernstein-type inequalities for $L^{3/2}$ observables of a non-mixing diffusion on the two-dimensional torus. The proof combines perturbation estimates for the Feynman--Kac resolvent with the stationary increment property of additive functionals of stationary processes.
\end{abstract} \noindent
 \vskip 2cm

\section{Introduction and Main Results}

Let $(X_t)_{t\ge0}$ be a jointly measurable conservative Markov process on a Polish space $E$, with invariant probability measure $\mu$ and transition semigroup $(P_t)_{t\ge0}$. We are concerned with deviations of the time averages
\[
\frac1t\int_0^t f(X_r) \,\d r
\]
from $\mu(f):=\int_E f\,\d\mu$. For an initial distribution $\nu$, we write $\mathbb P^\nu$ and $\mathbb E^\nu$ for the corresponding probability and expectation. Our main interest is in Bernstein-type bounds of the form
\begin{equation}\label{eq:intro:Bernstein}
\mathbb P^\mu\left(
\left|\frac1t\int_0^t f(X_r)\,\d r\right|\ge\delta
\right)
\le C\exp\left(
-\frac{t\delta^2}{a(f)+b(f)\delta}
\right),
\qquad \mu(f)=0,
\end{equation}
where the quantities $a(f)$ and $b(f)$ reflect both the dynamics of the process and the size of $f$. Such inequalities are also used to estimate convergence rates of empirical measures in Wasserstein distance;
see \cite{SWWZ,TWZ,WWZ}. In this paper, we aim to find conditions on the generator that yield such bounds and allow $a(f)$ and $b(f)$ to be controlled by suitable norms of $f$.

Exponential deviation bounds for additive functionals of Markov processes have been obtained through spectral
and functional analytic methods. Lezaud \cite{L01} used perturbation theory for linear
operators to establish Chernoff and Berry--Ess\'een inequalities under a spectral gap assumption.
Wu \cite{Wu} obtained a deviation inequality for non-reversible Markov processes by controlling the Feynman--Kac semigroup through the symmetrized Dirichlet
form. This connection with functional inequalities was further developed by Gao, Guillin and Wu \cite{GGW}, who established
Bernstein-type bounds for unbounded observables of symmetric Markov processes. Wu's inequality was also used in \cite{WWZ} to derive
Bernstein-type bounds for a class of diffusion processes.

For non-reversible processes, the symmetrized Dirichlet form need not capture the contribution of the antisymmetric part of the generator.
One way to retain this contribution is to use hypocoercivity in a modified inner product. Birrell and Rey-Bellet \cite{BR} and Monmarch\'e
\cite{M22} combined this approach with Feynman--Kac estimates to obtain Bernstein-type inequalities for hypocoercive Markov processes, including underdamped
Langevin dynamics under suitable assumptions. These results apply beyond the setting of a Poincar\'e inequality for the original symmetrized Dirichlet form,
but still rely on quantitative convergence to equilibrium in a modified $L^2$ norm.

A weaker condition is the iterated Poincar\'e inequality. To describe it, write $L^p(\mu):=L^p(E,\mu)$ equipped with norm $\|\cdot\|_p$, and assume that $(P_t)_{t\ge0}$ is strongly continuous on $L^2(\mu)$, with generator $(\L,\mathscr D(\L))$. Let $L^2_0(\mu)$ and $\mathscr D_0(\L)$ denote the
mean-zero subspaces of $L^2(\mu)$ and $\mathscr D(\L)$, respectively. The iterated Poincar\'e inequality requires that, for
some $\eta>0$,
\begin{equation}\label{eq:intro:IP}
\|\L u\|_2\ge\eta\|u\|_2,
\qquad u\in\mathscr D_0(\L).
\end{equation}
Unlike the usual Poincar\'e inequality, this condition controls $\L u$ rather than the quadratic form
$\langle-\L u,u\rangle_2$. In particular, it may hold even when the semigroup does not converge to equilibrium.

By the classical Bernstein form, we mean a bound of the form \eqref{eq:intro:Bernstein} with
\[
a(f)=c_1\|f\|_2^2,
\qquad
b(f)=c_2\|f\|_\infty,
\]
where $c_1,c_2>0$ are independent of $f$, $t$ and $\delta$. Huang and Li \cite{HL} established Bernstein-type inequalities under an iterated Poincar\'e inequality, first for discrete-time Markov chains and then for continuous-time Markov processes. Their bounds involve $\|f\|_2\|f\|_\infty$ in place of
$\|f\|_2^2$ in the Gaussian term. They sought to obtain the classical Bernstein form under the same assumption. Our main result gives this form up to numerical constants. More precisely, \eqref{eq:intro:IP} implies that, for every bounded mean-zero $f$, $t>0$ and $\delta>0$,
\begin{equation}\label{eq:intro:IPtail}
\mathbb P^\mu\left(
\left|\frac1t\int_0^t f(X_r)\,\d r\right|\ge\delta
\right)
\le
4\exp\left(
-\frac{\eta t\delta^2}
{64\|f\|_2^2+16\|f\|_\infty\delta}
\right),
\end{equation}
without assuming exponential convergence to equilibrium. We also prove the corresponding discrete-time result under $\|(I-P)u\|_2\ge\eta\|u\|_2$.
For other Bernstein inequalities for discrete-time Markov chains, see e.g. \cite{JSF}. This improvement also strengthens the results of \cite{SWWZ},
where Huang and Li's inequality was used under assumption \textbf{(H2)}. Our bound \eqref{eq:intro:IPtail} allows \textbf{(H3)}
to be replaced by \textbf{(H2)} in convergence rates proved therein, in particular for suitable underdamped Langevin dynamics.

The starting point of our method is a uniform bound on the resolvent
\[
R_s:=(s-\L)^{-1}
=\int_0^\infty \e^{-st}P_t\,\d t,
\qquad s>0,
\]
acting on $L^2_0(\mu)$. The iterated Poincar\'e inequality is equivalent,
up to constants, to
\[
\sup_{s>0}\|R_s\|_{L^2_0(\mu)\to L^2_0(\mu)}<\infty.
\]
This observation suggests replacing the $L^2$ resolvent bound by an estimate between spaces adapted to the observable. We introduce the following assumption involving a Gelfand triple.
\begin{asm}[Gelfand-triple resolvent bound] \label{AS0}
Let $\V$ be a Banach space with dual $\V^*$ such that
\[
\V \hookrightarrow L^2_0(\mu) \hookrightarrow \V^*
\]
continuously, with the first embedding dense. We denote the norms on $\V$ and $\V^*$ by $\|\cdot\|_{\V}$
and $\|\cdot\|_{\V^*}$, respectively. Suppose that, for every $s>0$, the resolvent $R_s$ admits a bounded extension $\widetilde R_s : \V^* \to \V$. More precisely, let $i : \V \to L^2_0(\mu)$ denote the embedding and define $j : L^2_0(\mu) \to \V^*$ by
\begin{equation} \label{def of j}
\< jf, g \>_{\V^*, \V}:=\< f,ig \>_2,
\qquad f\in L^2_0(\mu),\, g \in \V,
\end{equation}
where $\langle\cdot,\cdot\rangle_{\V^*, \V}$ denotes the dual pairing between $\V^*$ and $\V$. We require
\[
i \widetilde R_s j = R_s\quad\text{on }L^2_0(\mu),
\]
and
\begin{equation}
\sup_{s>0} \| \widetilde R_s \|_{\V^* \to \V} \le K.
\end{equation}
for some $K>0$.
\end{asm}

\begin{rem}\label{rem:IP''}
Identifying $\V$ with its image in $L^2_0(\mu)$, Assumption \ref{AS0} implies $\mathscr D_0(\L)\subset\V$ and
\begin{equation} \label{eq:IP''}
\|u\|_\V \le K\|j\L u\|_{\V^*},
\quad u\in\mathscr D_0(\L).
\end{equation}
Indeed, compatibility gives
$u=R_s(s-\L)u=i\widetilde R_sj(s-\L)u\in i(\V)$.
Thus
\[
\|u\|_\V
\le K\bigl(s\|ju\|_{\V^*}+\|j\L u\|_{\V^*}\bigr).
\]
Letting $s\downarrow0$, we obtain \eqref{eq:IP''}. In particular, by the continuous embeddings,
\[
\|u\|_2 \le K\|i\|_{\V \to L^2_0(\mu)} \|j\|_{L^2_0(\mu) \to \V^*} \cdot \|\L u\|_2,
\qquad u\in\mathscr D_0(\L).
\]
Thus Assumption \ref{AS0} implies the so-called iterated Poincar\'e inequality. When $\V=L^2_0(\mu)$, the converse also holds, up to constants;
see Lemma~\ref{lem:IP} below.
\end{rem}

For a bounded measurable function $f: E \to \R$ such that $\mu(f) = 0$, define
\[
\alpha_\V(f) := \sup_{\| u \|_\V\le 1} |\mu(fu)|, \quad \beta_\V(f) := \sup_{\|u\|_\V \le 1, \,  \| v \|_\V \le 1} |\mu(fuv)|.
\]
Both quantities are finite by the continuous embedding $\V \hookrightarrow L^2_0(\mu)$. If $\V$ is a Hilbert space, real polarization gives
\begin{equation} \label{beta}
\beta_\V(f)=\sup_{\|u\|_\V\le1}|\mu(fu^2)|.
\end{equation}
These quantities admit the following interpretations. Recalling \eqref{def of j}, we have
\begin{equation} \label{alpha}
\alpha_\V(f) = \| j f \|_{\V^*}.
\end{equation}
Let $Q$ denote the orthogonal projection onto $L^2_0(\mu)$,
given by
\begin{equation} \label{def of Q}
Qg:=g-\mu(g),\qquad g \in L^2(\mu).
\end{equation}
Let $M_f$ denote the multiplication operator on $L^2(\mu)$ given by $M_fg:=fg$. Define
$M_f^0:L^2_0(\mu)\to L^2_0(\mu)$ and
$\widehat M_f^0:\V\to \V^*$ by
\[
M_f^0:=QM_f\big|_{L^2_0(\mu)},
\qquad
\widehat M_f^0:=j\, M_f^0 \,i.
\]
Then
\[
\< \widehat M_f^0u,v \>_{\V^*, \V}=\mu(fuv),
\qquad u,v\in \V,
\]
and hence
\begin{equation} \label{def of b}
\beta_\V(f) = \| \widehat M_f^0 \|_{\V \to \V^*}.
\end{equation}
Now we can state our main result. Throughout the paper, we exclude the trivial case where $f=0$ $\mu$-almost everywhere.

\begin{thm} \label{main}
With the notation introduced above, under Assumption \ref{AS0}, let $f: E \to \R$ be a bounded measurable function such that $\mu(f) = 0$. Then for any $t>0$ and $\theta \in \R$ satisfying $4K\beta_\V(f)|\theta|\le1$,
\begin{equation} \label{eq:moment}
\E^\mu \left[ \e^{\theta \int_0^{t} f(X_r) \, \d r}  \right] \le 2 \e^{8K \alpha_\V(f)^2 \theta^2 t}.
\end{equation}
Consequently, for any $\delta > 0$,
\begin{equation} \label{eq:tail}
\mathbb P^\mu \left( \left|\frac1t \int_0^{t} f(X_r) \, \d r  \right| \ge \delta  \right) \le 4 \exp \left( - \frac{t\delta^2}{8K \big[ 4\alpha_\V(f)^2 + \beta_\V(f) \delta \big] } \right).
\end{equation}
More generally, if $\nu \ll \mu$ and $\frac{\d \nu}{\d \mu} \in L^p(\mu)$ for some $1 < p \le \infty$, then, with $q$ the H\"older conjugate of $p$,
\begin{equation} \label{eq:tail2}
\mathbb P^\nu \left( \left|\frac1t \int_0^{t} f(X_r) \, \d r  \right| \ge \delta  \right) \le 4^{\frac{1}{q}} \left\| \frac{\d \nu}{\d \mu} \right\|_p \exp \left( - \frac{t\delta^2}{8Kq \big[ 4\alpha_\V(f)^2 + \beta_\V(f) \delta \big]} \right).
\end{equation}
\end{thm}

We prove this result in Section~\ref{sec:main}. Under stationarity, additive functionals have stationary increments. We use this property to reduce the exponential moment estimate to an integral bound, which we prove by perturbation estimates for the Feynman--Kac resolvent. In Section~\ref{section:fin}, we first apply Theorem~\ref{main} under an iterated Poincar\'e inequality and also prove a discrete-time analogue. We then consider Poincar\'e and Sobolev--Poincar\'e inequalities and recover the Bernstein-type bounds in \cite{GGW,WWZ}, up to numerical constants. These latter applications illustrate how the resolvent approach also yields known bounds. In Section~\ref{sec:non-mixing}, we give a criterion for our resolvent assumption in terms of a positive self-adjoint operator $A$, with $\V=\mathscr D(A^{1/2})$ and $\|u\|_\V=\|A^{1/2}u\|_2$. For a normal generator satisfying the iterated Poincar\'e inequality, we may take $A=|\L|=(\L^*\L)^{1/2}$ on $L^2_0(\mu)$.
This choice captures the contribution of the antisymmetric part of the generator to the resolvent estimate. We illustrate this with a simple but instructive example:
the non-mixing diffusion on $\mathbb T^2$ generated by
\[
\L=\partial_x+\partial_y^2.
\]
Although its symmetrized Dirichlet form vanishes on functions depending only on the first coordinate, our resolvent estimate yields a Bernstein-type inequality
for every mean-zero $f\in L^{3/2}(\mu)$.

\section{Proof of Theorem \ref{main}} \label{sec:main}

We first show that, to prove Theorem \ref{main}, it suffices to establish the exponential moment bound \eqref{eq:moment}. Indeed, assume that \eqref{eq:moment} holds. For $\delta > 0$, choose
\[
\theta
:=\frac{\delta}
{4K\big[4\alpha_\V(f)^2+\beta_\V(f)\delta \big]}.
\]
Then $4 K\beta_\V(f)\theta\le 1$ and $8K\alpha_\V(f)^2\theta\le\delta/2$. By the exponential Markov inequality and \eqref{eq:moment},
\begin{align*}
\mathbb P^\mu\left( \frac1t\int_0^t f(X_r)\,\d r\ge\delta \right)
&\le
\e^{-\theta t\delta}
\E^\mu\left[
\e^{\theta\int_0^t f(X_r)\,\d r}
\right] \le
2\exp\big(-t\theta\delta+8K\alpha_\V(f)^2\theta^2t\big)\\
&\le
2\exp\left(
-\frac{t\delta^2}
{8K\big[4\alpha_\V(f)^2+\beta_\V(f)\delta \big]}
\right).
\end{align*}
where the last inequality follows from
\[
-t\theta\delta+8K\alpha_\V(f)^2\theta^2t = t\theta \big( -\delta + 8K\alpha_\V(f)^2\theta \big) \le -\frac{t\theta\delta}{2}.
\]
Since $\alpha_\V(-f)=\alpha_\V(f)$ and $\beta_\V(-f)=\beta_\V(f)$, applying the same argument to $-f$ and adding the two tail bounds proves \eqref{eq:tail}. Finally, let $\nu\ll\mu$ with $\frac{\d \nu}{\d \mu} \in L^p(\mu)$, where $1<p\le\infty$, and let $q$ be the H\"older conjugate of $p$. By H\"older's inequality,
\[
\mathbb P^\nu \left( \left|\frac1t \int_0^{t} f(X_r) \, \d r  \right| \ge \delta  \right) \le \left\| \frac{\d \nu}{\d \mu} \right\|_p \cdot \mathbb P^\mu \left( \left|\frac1t \int_0^{t} f(X_r) \, \d r  \right| \ge \delta  \right)^{\frac1q}.
\]
Thus \eqref{eq:tail2} follows from \eqref{eq:tail}. It remains to prove \eqref{eq:moment}. The following lemma reduces this task to estimating
Laplace integrals of exponential moments. It requires neither the Markov property nor independent increments.

\begin{lem} \label{lem:pointwise}
Let $(S_t)_{t\ge 0}$ be a jointly measurable real-valued stochastic process with $S_0=0$ and stationary increments, i.e., $S_{t+r}-S_r$ has the same distribution as $S_t$ for all $t,r\ge0$. For $\theta \in \R$, set $F_\theta(t) := \E [\e^{\theta S_t}]$. If, for some $z > 0$ and $C_+, C_-  > 0$,
\[
z \int_{0}^{\infty} \e^{-zs} F_{2\theta} (s) \, \d s \le C_+, \quad z \int_{0}^{\infty} \e^{-zs} F_{-2\theta} (s) \, \d s \le C_-,
\]
then for every $t \ge 0$,
\begin{equation} \label{eq:pointwise}
F_\theta(t) \le \sqrt{C_+ C_-} \, \e^{\frac{zt}{2}}.
\end{equation}
\end{lem}

\begin{proof}
By stationary increments and the Cauchy-Schwarz inequality, for any $t, s \ge 0$,
\[
F_\theta(t) = \E \left[  \e^{\theta(S_{t+s} - S_s)}  \right] \le \sqrt{\E \big[ \e^{2\theta S_{t+s}} \big] \E \big[ \e^{-2\theta S_s } \big]} = \sqrt{F_{2\theta}(t+s) F_{-2\theta}(s)}.
\]
Integrating this inequality against $z \e^{-zs} \, \d s$ and applying the Cauchy-Schwarz inequality again, we obtain
\begin{align*}
F_\theta(t) &\le \left( z \int_{0}^{\infty} \e^{-zs} F_{2\theta} (t+s) \, \d s  \right)^{1/2} \cdot \left( z \int_{0}^{\infty} \e^{-zs} F_{-2\theta} (s) \, \d s  \right)^{1/2}\\
& \le \sqrt{C_-}  \left( z \int_{0}^{\infty} \e^{-zs} F_{2\theta} (t+s) \, \d s  \right)^{1/2}.
\end{align*}
Since
\[
z \int_{0}^{\infty} \e^{-zs} F_{2\theta} (t+s) \, \d s = \e^{z t} \cdot z \int_{t}^{\infty} \e^{-zr} F_{2\theta} (r) \, \d r \le C_+ \e^{z t},
\]
the desired bound follows.
\end{proof}
Under $\mathbb P^\mu$, the process $(X_t)_{t\ge0}$ is stationary, so the additive functional
\[
S_t:=\int_0^t f(X_r)\,\d r
\]
has stationary increments. For $\theta \in \R$, define
\[
F_\theta(t) : = \E^\mu \big[ \e^{\theta \int_0^t f(X_r)\,\d r}  \big], \qquad t \ge 0.
\]
By Lemma \ref{lem:pointwise}, it therefore suffices to establish the corresponding Laplace integral bounds for $F_{2\theta}$ and $F_{-2\theta}$.

\subsection{Resolvent estimates on $L^2(\mu)$}

The exponential moment in the theorem is associated with the Feynman--Kac semigroup generated by $\mathcal L + \theta M_f$. To estimate its resolvent, we separate the constant and mean-zero components in $L^2(\mu)$. This leads us to first consider $s-\mathcal L-\theta M_f^0$ on $L^2_0(\mu)$. Fix $s > 0$ and $\theta \in \R$ such that $K\beta_\V(f)|\theta|<1$. By Assumption \ref{AS0} and \eqref{def of b},
$$\|\widetilde R_s \widehat M_f^0 \|_{\V \to \V} \le K \beta_\V(f). $$
Thus the Neumann series defines a bounded operator $R_{s, \theta} : \V^*\to \V$ given by
\[
R_{s, \theta} := (I-\theta \widetilde R_s \widehat M_f^0)^{-1} \widetilde R_s
\]
with
\begin{equation} \label{bounds of r}
\| R_{s, \theta} \|_{\V^* \to \V} \le \frac{K}{1 -  K \beta_\V(f)|\theta|}.
\end{equation}
This operator agrees with $(s - \L - \theta M_f^0)^{-1}$ on $L^2_0(\mu)$, namely
\[
i R_{s, \theta} j = (s - \L - \theta M_f^0)^{-1}.
\]
Indeed, given any $g \in L_0^2(\mu)$, we write $u = R_{s, \theta} j g \in \V$. By definition,
\[
u - \theta \widetilde R_s \widehat M_f^0 u = \widetilde R_s j g.
\]
Since $\widehat M_f^0 = j M_f^0 i$,
\[
u = \widetilde R_s j (g + \theta M_f^0 i u).
\]
By the compatibility $i \widetilde R_s j = R_s$, we have
\[
i u = R_s  (g + \theta M_f^0 i u).
\]
Therefore, $i u \in \mathscr D_0(\L)$, and
\[
(s - \L - \theta M_f^0) i u = g.
\]
This shows that $s-\L-\theta M_f^0$ maps $\mathscr D_0(\L)$ onto $L^2_0(\mu)$. To prove injectivity, suppose that $v\in\mathscr D_0(\L)$ satisfies
\[
(s-\L-\theta M_f^0)v=0.
\]
Then
\[
v=\theta R_s M_f^0v
  =i\bigl(\theta\widetilde R_s j M_f^0v\bigr),
\]
so $v=iw$ for some $w\in \V$. Since $i$ is injective,
\[
w=\theta\widetilde R_s\widehat M_f^0w,
\]
The bound $K\beta_V(f)|\theta|<1$ implies that $w=0$, and hence $v=0$. Thus
\begin{equation} \label{res}
(s-\L-\theta M_f^0)^{-1}=iR_{s,\theta}j
\quad\text{on }L^2_0(\mu).
\end{equation}

Now we define the function
\[
\chi_f(s, \theta) := \<j f, R_{s, \theta} j f  \>_{\V^*, \V} = \mu(f \, i R_{s, \theta} j f).
\]
Then by \eqref{alpha} and \eqref{bounds of r},
\begin{equation} \label{eq:chi}
| \chi_f(s, \theta) | \le  \| R_{s, \theta}  \|_{\V^* \to \V} \| j f \|_{\V^*}^2 \le \frac{K \alpha_\V(f)^2}{1 -  K \beta_\V(f)|\theta|}.
\end{equation}
We henceforth suppress the canonical embeddings $i$ and $j$, and write $\chi_f(s, \theta) = \mu(f R_{s, \theta} f)$. Using the decomposition $L^2(\mu) = \R \mathbf 1 \oplus L^2_0(\mu)$ where $\mathbf 1$ denotes the constant function equal to one,  we seek $c\in\R$ and $u\in\mathscr D_0(\L)$ such that
\[
(s-\L-\theta M_f)(c+u)=d+g,
\]
for given $d\in\R$ and $g\in L^2_0(\mu)$. Taking the mean-zero part and the mean, and using \eqref{res}, we obtain
\[
u= R_{s,\theta}(g+\theta cf),
\qquad
sc-\theta\mu(fu)=d.
\]
If $s - \theta^2 \chi_f(s, \theta) \ne 0$, these equations give
\[
c = \frac{d + \theta \mu(f R_{s,\theta}g)}{s - \theta^2 \chi_f(s, \theta)}, \qquad u= R_{s,\theta}(g+\theta cf).
\]
Since $iR_{s,\theta}j$ is bounded on $L^2_0(\mu)$ and $s-\theta^2\chi_f(s,\theta)\ne 0$, these formulas define a bounded inverse of $s-\L-\theta M_f$ on the full $L^2(\mu)$, with values in $\D(\L)$. In particular, taking $d= 1$ and $g = 0$, we obtain
\[
\mu\big( (s-\L-\theta M_f)^{-1} \mathbf 1  \big) = \frac{1}{s - \theta^2 \chi_f(s, \theta)}.
\]
The density of $\V$ in $L^2_0(\mu)$ implies that $\alpha_\V(f) > 0$. Suppose that $\theta \ne 0$ and $K \beta_\V(f) |\theta| \le 1/2$. Then by \eqref{eq:chi}, for every $s \ge 4K \alpha_\V(f)^2 \theta^2$,
\begin{equation} \label{control of w}
0 < \mu\big( (s-\L-\theta M_f)^{-1} \mathbf 1  \big) \le \frac{2}{s}.
\end{equation}

\subsection{From resolvent estimates to exponential moment bounds}

From now on, we write $G_\theta:=\L+\theta M_f$ with domain $\D(\L)$. Let $(T_t^\theta)_{t\ge0}$ be the Feynman--Kac semigroup
\[
T_t^\theta h(x)
:= \E^x\left[ \e^{\theta\int_0^t f(X_r)\,\d r} h(X_t) \right].
\]
Since $f$ is bounded, this is a positive strongly continuous semigroup on $L^2(\mu)$ with generator
$G_\theta$. Moreover, it follows from $| T_t^\theta h(x)| \le \e^{ \| f \|_\infty |\theta| t} P_t|h| (x)$ that
\[
\|T_t^\theta\|_{L^2(\mu)\to L^2(\mu)}
\le \e^{\|f\|_\infty |\theta| t}.
\]
Thus, for $s>\|f\|_\infty |\theta|$, the Laplace representation
\[
(s-G_\theta)^{-1}
=\int_0^\infty \e^{-st}T_t^\theta\,\d t
\]
shows that $(s-G_\theta)^{-1}$ is positive. The preceding subsection gives invertibility for every $s\ge 4K\alpha_\V(f)^2\theta^2$, provided that
$\theta\ne0$ and $K\beta_\V(f)|\theta|\le1/2$. However, the operator norm bound above does not justify the Laplace representation throughout this range. We will need positivity of the resolvent to discard a nonnegative remainder in the Laplace integral estimate. Here an operator is called positive if it preserves
nonnegativity. The following lemma extends positivity to all $s\ge4K\alpha_\V(f)^2\theta^2$.

\begin{lem}[Continuation of positivity]\label{lem:positive}
Let $G$ be the generator of a positive $C_0$ semigroup $(T_t)_{t\ge0}$ on $L^2(\mu)$. Suppose that, for some $s_0\in\R$, the operator $s-G$
has a bounded inverse on $L^2(\mu)$ for every $s\ge s_0$. Then $(s-G)^{-1}$ is positive for every $s\ge s_0$.
\end{lem}

\begin{proof}
Choose $M \ge 1$ and $\omega \in \R$ such that $\| T_t \|_{L^2(\mu) \to L^2(\mu)} \le M \e^{\omega t}$ for all $t\ge 0$. Fix $s_1 > \max \{s_0, \omega \}$. For every $s\ge s_1$,
the Laplace representation
\[
(s-G)^{-1}=\int_0^\infty \e^{-st}T_t\,\d t
\]
shows that $(s-G)^{-1}$ is positive. Norm continuity of the resolvent gives
\[
C_1 := \sup_{s \in [s_0, s_1]} \| (s - G)^{-1} \|_{L^2(\mu) \to L^2(\mu)} < \infty.
\]
If $(s-G)^{-1}$ is positive for some $s\in[s_0,s_1]$, then for $\delta \in (0, C_1^{-1})$ with $s-\delta\ge s_0$, the resolvent identity
\[
(s - \delta -G)^{-1} = \sum_{j=0}^\infty\delta^j(s-G)^{-j-1}
\]
yields that $(s - \delta -G)^{-1}$ is positive. Iterating this argument from $s_1$ down to $s_0$, with steps smaller than $C_1^{-1}$, proves the assertion.
\end{proof}

Suppose that $\theta\ne0$ and $K\beta_\V(f)|\theta|\le1/2$, and set $s_0:=4K\alpha_\V(f)^2\theta^2$. Applying Lemma~\ref{lem:positive} to $G_\theta$, we deduce
that $(s-G_\theta)^{-1}$ is positive for every $s\ge s_0$. Set
\[
w := (s_0 - G_\theta)^{-1} \mathbf 1.
\]
Then $w\ge 0$, and \eqref{control of w} gives
\begin{equation} \label{integ w}
\mu(w) \le \frac{2}{s_0}.
\end{equation}
Since $w\in \D(G_\theta)$ and $(s_0 - G_\theta) w = \mathbf 1$, we have
\[
\frac{\d}{\d t} \big(\e^{-s_0 t}T_t^\theta w \big) = -\e^{-s_0 t} T_t^\theta \mathbf 1.
\]
Therefore, for any $\mathrm T > 0$,
\begin{equation} \label{integ}
\int_{0}^{\mathrm T} \e^{-s_0 t} T_t^\theta \mathbf 1 \, \d t = w - \e^{-s_0 \mathrm T}T_{\mathrm T}^\theta w \le w,
\end{equation}
where the inequality follows from the positivity of $T_{\mathrm T}^\theta$ and $w\ge0$. For $t \ge 0$, recalling the definition of $F_\theta$,
\[
F_\theta (t)  = \E^\mu \big[ \e^{\theta\int_0^t f(X_r)\,\d r} \big] = \mu(T_t^\theta \mathbf 1).
\]
Then integrating \eqref{integ} against $\mu$ yields
\[
\int_{0}^{\mathrm T} \e^{-s_0 t} F_\theta (t) \, \d t  \le \mu(w).
\]
Letting $\mathrm T \to \infty$ and using \eqref{integ w} and monotone convergence, we obtain
\begin{equation}
s_0 \int_{0}^{\infty} \e^{-s_0 t} F_\theta (t) \, \d t \le 2.
\end{equation}

Now suppose that $\theta\ne0$ and $4K\beta_\V(f)|\theta|\le1$. Applying the preceding bound with $2\theta$ and $-2\theta$ in place of $\theta$, and setting
$z = 16K\alpha_\V(f)^2\theta^2$, we obtain
\[
z\int_0^\infty \e^{-zt}F_{2\theta}(t)\,\d t\le2,
\qquad
z\int_0^\infty \e^{-zt}F_{-2\theta}(t)\,\d t\le2.
\]
Lemma~\ref{lem:pointwise} therefore gives
\[
F_\theta(t)\le2\e^{zt/2} =2\e^{8K\alpha_\V(f)^2\theta^2t}, \qquad t \ge 0,
\]
which proves \eqref{eq:moment}. The case $\theta=0$ is immediate. The proof is complete.

\subsection{Extension to unbounded observables}

For later applications, we extend Theorem~\ref{main} to unbounded observables by approximation.

\begin{lem} \label{lem:approx}
Under Assumption \ref{AS0}, let $f \in L^1(\mu)$ with $\mu(f) = 0$. Suppose there exist bounded measurable functions $(f_n)_{n \ge 1}$ with $\mu(f_n) = 0$ and constants $\alpha^*, \beta^* \in [0, \infty)$ such that
\[
\| f_n - f \|_1 \to 0, \qquad \limsup_{n \to \infty} \alpha_{\V}(f_n) \le \alpha^*, \qquad \limsup_{n \to \infty} \beta_{\V}(f_n) \le \beta^*.
\]
Then for any $t>0$ and $\theta \in \R$ satisfying $4K \beta^* |\theta| \le 1$,
\begin{equation} \label{eq:moment''}
\E^\mu \left[ \e^{\theta \int_0^{t} f(X_r) \, \d r}  \right] \le 2 \e^{8K (\alpha^*)^2 \theta^2 t}.
\end{equation}
The tail bounds \eqref{eq:tail} and \eqref{eq:tail2} also hold, with $\alpha_\V(f)$ and $\beta_\V(f)$ replaced by $\alpha^*$ and $\beta^*$, respectively.
\end{lem}

\begin{proof}
As before, it suffices to prove the exponential moment bound \eqref{eq:moment''}. By stationarity,
\[
\E^\mu \left[ \int_0^t |f(X_r)| \, \d r  \right] = t \| f \|_1 < \infty, \quad \E^\mu \left[ \bigg| \int_0^t f_n(X_r) \, \d r - \int_0^t f(X_r) \, \d r \bigg| \right] \le t \| f_n - f \|_1 \to 0.
\]
Thus for each fixed $t$, $\int_0^t f(X_r) \, \d r$ is finite almost surely, and a subsequence $\bigl(\int_0^t f_{n_k}(X_r)\,\d r\bigr)_{k \ge1}$ converges to it
almost surely. If $4K \beta^* |\theta| < 1$, then by our assumption, $4K  \beta_{\V}(f_n) |\theta| \le 1$ for all sufficiently large $n$. Along that subsequence, Fatou's Lemma and \eqref{eq:moment} imply
\[
\E^\mu \left[ \e^{\theta \int_0^{t} f(X_r) \, \d r}  \right] \le \liminf_{k \to \infty} \E^\mu \left[ \e^{\theta \int_0^{t} f_{n_k}(X_r) \, \d r}  \right] \le 2 \e^{8K (\alpha^*)^2 \theta^2 t}.
\]
The endpoint case $4K\beta^*|\theta|=1$ follows by applying again Fatou's lemma with $\theta_m=(1-m^{-1})\theta$, $m\ge2$. The proof is complete.
\end{proof}

\begin{rem} \label{rem:sobolev}
Suppose in addition that, for some $r\in(2,\infty]$, $\V \hookrightarrow L^r(\mu)$, i.e., there exists $C_r > 0$ such that
\begin{equation} \label{eq:s}
\|u\|_r\le C_r\|u\|_\V,\qquad u\in\V.
\end{equation}
Then Lemma~\ref{lem:approx} allows us to extend the Bernstein bounds to mean-zero observables in $L^p(\mu)$, where $p$ is the H\"older conjugate of $r/2$.
Indeed, let $f\in L^p(\mu)$ with $\mu(f)=0$. We use the same definitions of $\alpha_\V(f)$ and
$\beta_\V(f)$; the following estimates ensure their finiteness. By H\"older's inequality and \eqref{eq:s}, for every $u, v\in\V$,
\[
|\mu(fu)| \le \| f \|_{r'} \| u \|_{r} \le C_r \| f \|_p \| u \|_\V, \ \ |\mu(f u v)| \le \| f \|_p \| u \|_{r} \| v \|_{r}\le C_r^2 \| f \|_p \| u \|_\V \| v \|_\V.
\]
Therefore, by the definition of $\alpha_\V$ and $\beta_\V$,
\begin{equation} \label{eq:control'}
\alpha_\V(f) \le C_r \| f \|_p, \qquad \beta_\V(f) \le C_r^2 \| f \|_p.
\end{equation}
Consider the mean-zero truncations
\[
f_n := (-n) \vee (f \wedge n) - \mu \big( (-n) \vee (f \wedge n)   \big), \qquad n \ge 1.
\]
Then the dominated convergence theorem gives
\begin{equation} \label{eq:approx'}
\| f - f_n \|_p \le 2 \big\| (-n) \vee (f \wedge n) - f \big\|_p \to 0 \ \ \text{as} \ \ n \to \infty.
\end{equation}
Applying \eqref{eq:control'} to $f_n - f$ and using
\[
|\alpha_\V(f_n)-\alpha_\V(f)| \le\alpha_\V(f_n-f), \qquad |\beta_\V(f_n)-\beta_\V(f)| \le\beta_\V(f_n-f),
\]
we obtain
\begin{equation*}
\alpha_\V(f_n) \to \alpha_\V(f), \qquad \beta_\V(f_n) \to \beta_\V(f) \ \ \text{as} \ \ n \to \infty.
\end{equation*}
Since \eqref{eq:approx'} also implies convergence in $L^1(\mu)$, Lemma~\ref{lem:approx} applies with $\alpha^*=\alpha_\V(f)$ and $\beta^*=\beta_\V(f)$.
\end{rem}

\section{Bernstein-type inequalities under functional inequalities} \label{section:fin}

\subsection{Iterated Poincar\'e inequalities}

Bernstein-type inequalities for Markov chains under iterated Poincar\'e inequalities are studied in \cite{HL}. We say $(P_t)_{t \ge 0}$ satisfies the iterated Poincar\'e inequality (or admits an IP gap), if there exists a constant $\eta > 0$ such that
\begin{equation}\label{eq:IP}
\| \L u \|_2 \ge \eta \| u \|_2, \qquad u \in \D_0(\L).
\end{equation}
As noted in Remark~\ref{rem:IP''}, a uniform resolvent bound on $L^2_0(\mu)$ implies \eqref{eq:IP}. The following lemma proves the converse.

\begin{lem} \label{lem:IP}
Assume the iterated Poincar\'e inequality \eqref{eq:IP} holds. Recall that $R_s = (s - \L)^{-1}$. Then for any $s > 0$,
\[
\|R_s\|_{L^2_0(\mu)\to L^2_0(\mu)}\le s^{-1}.
\]
And for any $0 < s < \eta$,
\[
\|R_s\|_{L^2_0(\mu)\to L^2_0(\mu)}\le (\eta - s)^{-1}.
\]
In particular,
\begin{equation} \label{eq:uniform}
\sup_{s > 0} \| R_s\|_{L^2_0(\mu)\to L^2_0(\mu)} \le 2 \eta^{-1}.
\end{equation}
Moreover, $-\L : \D_0(\L) \to L_0^2(\mu)$ is onto and has a bounded inverse of norm at most $\eta^{-1}$.
\end{lem}

\begin{proof}
The first bound follows from the Laplace representation of $R_s$ and the contraction of $(P_t)_{t \ge 0}$. For $0 < s < \eta$, if $u = R_s g$, then
\[
\eta \| u \|_2 \le \| \L u \|_2 \le s \| u \|_2 + \| g \|_2,
\]
which proves the second bound. Combining these two bounds we immediately obtain \eqref{eq:uniform}.

The resolvent identity $R_s-R_r=(r-s)R_s R_r$ and \eqref{eq:uniform} imply that
\[
\|R_s-R_r\|_{L^2_0(\mu)\to L^2_0(\mu)}
\le 4\eta^{-2}|s-r|.
\]
Hence $R_s$ converges in operator norm to a bounded operator
$R_0$ as $s\downarrow0$. For every $g\in L^2_0(\mu)$, as $s\downarrow0$,
\[
R_sg \rightarrow R_0g,
\qquad
\L R_sg=sR_sg-g \rightarrow-g
\quad\text{in }L^2(\mu).
\]
Since $\L$ is closed, $R_0g\in\D_0(\L)$ and $-\L R_0g=g$. Thus $-\L$ is onto. Its injectivity follows directly from \eqref{eq:IP}, so
$R_0=(-\L)^{-1}$. Finally, applying \eqref{eq:IP} to $R_0g$ gives
\[
\eta\|R_0g\|_2\le\|\L R_0g\|_2=\|g\|_2,
\]
and hence $\|(-\L)^{-1}\|_{L^2_0(\mu)\to L^2_0(\mu)} \le\eta^{-1}$.
\end{proof}

Therefore when $\V=L^2_0(\mu)$, Assumption \ref{AS0} holds for $K = 2\eta^{-1}$ and we have
\[
\alpha_\V(f) = \| f \|_2, \qquad \beta_\V(f) \le \| f \|_\infty.
\]
Thus Theorem~\ref{main} gives the following consequence.

\begin{cor} \label{cor:BP}
Under the iterated Poincar\'e inequality \eqref{eq:IP}, let $f: E \to \R$ be a bounded measurable function such that $\mu(f) = 0$. Then for any  $t>0$ and $\delta > 0$,
\begin{equation} \label{eq:tail IP}
\mathbb P^\mu \left( \left|\frac1t \int_0^{t} f(X_r) \, \d r  \right| \ge \delta  \right) \le 4 \exp \left( - \frac{\eta t\delta^2}{64 \| f \|_2^2 + 16 \| f \|_\infty \delta } \right).
\end{equation}
The corresponding tail bound for an initial distribution $\nu\ll\mu$ with $\frac{\d\nu}{\d\mu}\in L^p(\mu)$, $1<p\le\infty$, follows from \eqref{eq:tail IP} by H\"older's inequality.
\end{cor}

The resolvent approach can also be adapted to discrete-time Markov chains, which are also considered in \cite{HL}. We now establish the corresponding Bernstein-type inequality under the discrete-time iterated Poincar\'e inequality. Let $P$ be a Markov transition operator with invariant measure $\mu$. In this context, we say  $P$ satisfies the iterated Poincar\'e inequality (or admits an IP gap), if there exists a constant $\eta > 0$ such that
\begin{equation}\label{eq:IP'}
\| (I - P) u \|_2 \ge \eta \| u \|_2, \qquad u \in L_0^2(\mu).
\end{equation}
Since $\| (I - P) u \|_2 \le 2 \|u\|_2$, we always have $\eta \le 2$. The same gap was used by Dinwoodie \cite{Dinwoodie}
to obtain deviation bounds for nonreversible Markov chains. His bounds uniform in the starting point involve an additional comparison condition on the transition kernel. We now state a Bernstein-type inequality assuming only the iterated Poincar\'e inequality.

\begin{cor}
Let $(X_n)_{n \ge 0}$ be a Markov chain on state space $E$ with transition operator $P$ and invariant probability measure $\mu$.  Under the iterated Poincar\'e inequality \eqref{eq:IP'}, let $f: E \to \R$ be a bounded measurable function such that $\mu(f) = 0$. Then for any  $n \in \N_+$ and $\theta \in \R$ satisfying $64 \eta^{-1} \| f \|_\infty |\theta| \le 1$,
\begin{equation} \label{eq:moment'}
\E^\mu \left[ \e^{\theta \sum_{j=0}^{n-1} f(X_j)}  \right] \le 2 \e^{40 \eta^{-1} \| f \|_2^2 \theta^2 n}.
\end{equation}
Consequently, for any $\delta > 0$,
\begin{equation} \label{eq:tail IP'}
\mathbb P^\mu \left( \left|\frac1n \sum_{j=0}^{n-1} f(X_j) \right| \ge \delta  \right) \le 4 \exp \left( - \frac{\eta n \delta^2}{160 \bigl[\|f\|_2^2 + \| f \|_\infty \delta \bigr]} \right).
\end{equation}
The corresponding tail bound for an initial distribution $\nu\ll\mu$ with $\frac{\d\nu}{\d\mu}\in L^p(\mu)$, $1<p\le\infty$, follows from \eqref{eq:tail IP'} by H\"older's inequality.
\end{cor}

\begin{proof}
The proof follows a discrete-time adaptation of the argument used to prove Theorem~\ref{main}. As before, it suffices to prove the exponential moment bound \eqref{eq:moment'}. The case $\theta=0$ is immediate. In what follows, we assume that $\theta\ne0$ and set $K_\theta:=\e^{\theta f}P$. Then
\begin{equation} \label{def of F}
F_\theta (n) := \E^\mu \big[ \e^{\theta \sum_{j=0}^{n-1} f(X_j)}  \big]=\mu(K_\theta^n \mathbf 1)
\end{equation}
Define $G = P -I$, $W_\theta = (\e^{\theta f} - 1) P$ and $W_\theta^0 = Q W_\theta |_{L^2_0(\mu)} : L^2_0(\mu) \to L^2_0(\mu)$, recalling $Q$ is the orthogonal projection onto $L^2_0(\mu)$ defined by \eqref{def of Q}. Then $K_\theta = I + G + W_\theta$. If $\|f\|_\infty |\theta| \le \eta/32 \le 1/16$, then Taylor expansion, the contraction property of $P$, and $\mu(f) = 0$ give
\begin{equation} \label{eq:theta1}
\| W_\theta \|_{L^2(\mu) \to L^2(\mu)} \le 2 \| f \|_\infty |\theta|, \qquad |\mu(W_\theta \mathbf 1)|  = |\mu(\e^{\theta f} - 1)| \le  \|f\|_2^2 \theta^2.
\end{equation}
Using the inequality $|\e^y - 1| \le |y| \e^{|y|}$ for $y \in \R$, we also have
\begin{equation} \label{eq:theta2}
\| Q W_\theta \mathbf 1 \|_2 \le  \| W_\theta \mathbf 1 \|_2  \le \e^{\| f\|_\infty |\theta|}  \|f\|_2 |\theta|.
\end{equation}
Moreover, for any $u \in L_0^2(\mu)$,
\begin{equation} \label{eq:theta3}
|\mu(W_\theta u)| \le \e^{\| f\|_\infty |\theta|} |\theta|  \mu(|f| \, |P u|)  \le  \e^{\| f\|_\infty |\theta|}  \|f\|_2 |\theta| \cdot \| u \|_2.
\end{equation}
Then we aim to define the inverse of $s - G - W_\theta$ on the full $L^2(\mu)$. For given $d\in \R$ and $g \in L_0^2(\mu)$, we seek $c\in \R$ and $u \in L_0^2(\mu)$ such that
\[
(s - G - W_\theta)(c + u) = d + g.
\]
Taking the mean-zero part and the mean, we obtain
\begin{equation} \label{eq:origin}
(s - G - Q W_\theta) u = g + c \, Q W_\theta \mathbf 1, \qquad \big( s - \mu(W_\theta \mathbf 1) \big) c - \mu(W_\theta u) = d.
\end{equation}
Now for $s > 0$, define the perturbed resolvent on $L^2_0(\mu)$ by
\[
\mathcal R_{s,\theta} :=\left(s-G- W_\theta^0 \right)^{-1}.
\]
Since $G=P-I$ generates a strongly continuous Markov contraction semigroup and satisfies the iterated Poincar\'e inequality \eqref{eq:IP'} with constant $\eta$,
Lemma~\ref{lem:IP} gives
\[
\sup_{s>0}\|(s-G)^{-1}\|_{L^2_0(\mu)\to L^2_0(\mu)}
\le\frac2\eta.
\]
On the other hand, by \eqref{eq:theta1}, we have $\| W_\theta^0 \|_{L_0^2(\mu) \to L_0^2(\mu)} \le 2 \| f \|_\infty |\theta|$. Since
\[
\|(s-G)^{-1} W_\theta^0\|_{L^2_0(\mu)\to L^2_0(\mu)}
\le\frac{4\|f\|_\infty|\theta|}{\eta}\le\frac18,
\]
a Neumann series gives
\[
\mathcal R_{s,\theta}
=\bigl(I-(s-G)^{-1}W_\theta^0 \bigr)^{-1}(s-G)^{-1}
\]
and hence
\begin{equation}  \label{eq:bound}
\|\mathcal R_{s,\theta}\|_{L^2_0(\mu)\to L^2_0(\mu)}
\le\frac{2/\eta}{1-4\|f\|_\infty|\theta|/\eta}
\le\frac{16}{7\eta},
\qquad s>0.
\end{equation}
Returning to \eqref{eq:origin}, we obtain, whenever the denominator below is nonzero,
\begin{equation}\label{eq:origin'}
c=
\frac{d+\mu(W_\theta\mathcal R_{s,\theta}g)}
{s-\mu(W_\theta\mathbf1)
-\mu(W_\theta\mathcal R_{s,\theta}QW_\theta\mathbf1)},
\qquad
u=\mathcal R_{s,\theta}(g+c\,QW_\theta\mathbf1).
\end{equation}
Combining \eqref{eq:theta1}, \eqref{eq:theta2}, \eqref{eq:theta3}, \eqref{eq:bound} and using $\eta \le 2$, we obtain
\begin{align*}
\big| \mu(W_\theta \mathbf 1) + \mu(W_\theta \mathcal R_{s, \theta} Q W_\theta \mathbf 1)  \big| &\le \| f \|_2^2 \theta^2 + \frac{16}{7\eta} \e^{2 \|f\|_\infty |\theta|} \| f \|_2^2 \theta^2\\
& \le \big(2+\frac{16}{7} \e^{\frac18} \big)  \eta^{-1} \| f \|_2^2 \theta^2
 < 10 \eta^{-1} \| f \|_2^2 \theta^2.
\end{align*}
Set $s_0 := 20 \eta^{-1} \|f\|_2^2 \theta^2$. Then for $s \ge s_0$, the denominator in \eqref{eq:origin'} is at least $s-s_0/2\ge s/2>0$. Thus these formulas define a bounded inverse of  $s-G-W_\theta=(1+s)I-K_\theta$ on $L^2(\mu)$. In particular,
\begin{equation} \label{control of w'}
0 < \mu \big( (1 + s -K_\theta)^{-1} \mathbf 1  \big) \le \frac2s.
\end{equation}
Since $K_\theta = \e^{\theta f} P$ is a bounded positive operator, $(\e^{t(K_\theta-I)})_{t\ge0}$ is a positive $C_0$ semigroup. Therefore Lemma~\ref{lem:positive} gives that, for any $s \ge s_0$, $(1 + s -K_\theta)^{-1}$ is positive. Set
\[
 w := (1 + s_0 -K_\theta)^{-1} \mathbf 1.
\]
Then $w \ge 0$ and \eqref{control of w'} gives
\begin{equation*}
\mu(w) \le \frac{2}{s_0}.
\end{equation*}
Iterating $(1+ s_0)w=\mathbf1+K_\theta w$, we obtain
\[
w=\sum_{k=0}^{N}\frac{K_\theta^k \mathbf1}{(1+s_0)^{k+1}} +\frac{K_\theta^{N+1}w}{(1 +  s_0)^{N+1}},
\qquad N\ge0.
\]
Since $w\ge0$ and $K_\theta$ is positive, the remainder is nonnegative. Integrating against $\mu$ and letting
$N\to\infty$, recalling the definition of $F_\theta$ we obtain by monotone convergence
\begin{equation} \label{eq:geo}
s_0 \sum_{k=0}^{\infty} (1+s_0)^{-k-1} F_\theta(k) \le s_0 \mu(w)\le2.
\end{equation}

We now use a discrete-time analogue of the argument in Lemma~\ref{lem:pointwise}. Under $\mathbb P^\mu$, the process $(X_n)_{n \ge 0}$ is stationary. By Cauchy-Schwarz inequality, for every $n \in \N_+$ and $k \ge 0$,
\begin{equation} \label{eq:sc}
F_\theta (n) \le \sqrt{F_{2\theta} (n+k) F_{-2\theta} (k)}.
\end{equation}
For any $\theta$ satisfying $64 \eta^{-1} \| f \|_\infty |\theta| \le 1$, applying \eqref{eq:geo} at $2\theta$ and $-2\theta$ with $s_0 = 80 \eta^{-1} \|f\|_2^2 \theta^2$, we have
\[
s_0 \sum_{k=0}^{\infty} (1+s_0)^{-k-1} F_{2\theta}(k) \le2, \qquad s_0 \sum_{k=0}^{\infty} (1+s_0)^{-k-1} F_{-2\theta}(k) \le2.
\]
Averaging \eqref{eq:sc} with respect to the geometric probabilities $s_0(1+s_0)^{-k-1}$ and applying the Cauchy--Schwarz inequality, we obtain
\begin{align*}
& F_\theta(n) \le \bigg(s_0\sum_{k=0}^\infty
(1+s_0)^{-k-1}F_{2\theta}(n+k)\bigg)^{1/2}
\bigg(s_0\sum_{k=0}^\infty
(1+s_0)^{-k-1}F_{-2\theta}(k)\bigg)^{1/2}\\
&\le \sqrt{2}\,
\bigg(s_0\sum_{k=0}^\infty
(1+s_0)^{-k-1}F_{2\theta}(n+k)\bigg)^{1/2} = \sqrt{2}\,(1+s_0)^{n/2}
\bigg(s_0\sum_{j=n}^\infty
(1+s_0)^{-j-1}F_{2\theta}(j)\bigg)^{1/2}\\
&\le 2(1+s_0)^{n/2}
\le 2\e^{40\eta^{-1}\|f\|_2^2\theta^2 n},
\end{align*}
which completes the proof.
\end{proof}

\subsection{Poincar\'e and Sobolev inequalities}

In this subsection, we derive Bernstein-type inequalities under Poincar\'e and Sobolev inequalities by choosing a suitable Banach space $\V$. Define the symmetrized Dirichlet form by
\[
\Ee_{\rm{sym}} (f, g) := \frac12 [\<-\L f, g  \>_2 + \<-\L g, f  \>_2 ], \qquad f, g \in \D(\L).
\]
Assume that $(\Ee_{\mathrm{sym}},\D(\L))$ is closable, and denote its closure by $(\Ee_{\mathrm{sym}},\D(\Ee_{\mathrm{sym}}))$. Then $\D(\L) \subset \D(\Ee_{\mathrm{sym}})$ and
\begin{equation}
\<-\L u, u  \>_2 = \Ee_{\rm{sym}} (u, u), \quad u \in \D(\L).
\end{equation}
Write $\D_0(\Ee_{\mathrm{sym}})$ for the mean-zero subspace of $\D(\Ee_{\mathrm{sym}})$. Suppose that, for some constant $C_P>0$,
the following Poincar\'e inequality holds:
\begin{equation}\label{eq:P}
\|u\|_2^2\le C_P\Ee_{\mathrm{sym}}(u,u), \qquad u\in\D_0(\Ee_{\mathrm{sym}}).
\end{equation}
Under these assumptions, a natural choice of $\V$ ensures that Assumption~\ref{AS0} holds, as shown below.

\begin{prp}
With the notation introduced above, under \eqref{eq:P}, set
\[
\V = \D_0(\Ee_{\mathrm{sym}}), \qquad \| u \|_\V = \Ee_{\rm{sym}} (u, u)^{1/2}.
\]
Then $\V$ is a Hilbert space such that $\V \subset L_0^2(\mu)$ continuously and densely, and Assumption~\ref{AS0} holds with $K =1$.
\end{prp}

\begin{proof}
Since $\D(\Ee_{\mathrm{sym}})$ is complete under the norm $\sqrt{\| u \|_2^2 + \Ee_{\mathrm{sym}}(u,u)}$, the Poincar\'e inequality \eqref{eq:P} implies that $\V$ is complete under $\|\cdot\|_\V$ and
\[
\|u\|_2\le\sqrt{C_P}\|u\|_\V,
\qquad u\in\V.
\]
Since $\D_0(\L)$ is dense in $L_0^2(\mu)$ and contained in $\V$, the embedding is dense.

We next verify Assumption~\ref{AS0}. For $g \in L_0^2(\mu)$ and $s > 0$, let $u = R_s g \in \D_0(\L)$. Taking the $L^2(\mu)$-inner product of $(s - \L)u = g$ with $u$, we obtain
\[
s \| u \|_2^2 + \| u \|_\V^2 = \< g, u \>_2 = \<j g, u  \>_{\V^*, \V} \le \| j g \|_{\V^*} \| u \|_\V.
\]
Consequently,
\begin{equation} \label{eq:boundness}
\| R_s g \|_\V \le \| j g \|_{\V^*}, \qquad s >0.
\end{equation}
The range $j(L^2_0(\mu))$ is norm dense in $\V^*$. Indeed, let $v_g\in\V$ be the Riesz representative of $jg$, so that
\[
(v_g,u)_\V = \langle g,iu\rangle_2,
\qquad u\in\V,
\]
where $(\cdot, \cdot)_\V$ denotes the inner product in $\V$. If $u\in\V$ is orthogonal to every $v_g$, taking $g=iu$ gives
$\|iu\|_2^2=0$, hence $u=0$. Thus $\{v_g:g\in L^2_0(\mu)\}$ is dense in $\V$, and the claim follows from the Riesz isomorphism. Since $j$ is injective, \eqref{eq:boundness} defines a contraction $jg \mapsto R_s g$ on this dense subspace. Therefore its unique continuous extension $\widetilde R_s: \V^* \to \V$ satisfies
\[
\| \widetilde R_s \|_{\V^* \to \V} \le 1.
\]
This proves Assumption~\ref{AS0} with $K=1$. The identity $i\widetilde R_sj=R_s$ holds by construction.
\end{proof}

Throughout the remainder of this subsection, we take $\V=\D_0(\Ee_{\mathrm{sym}})$ with the norm defined above. Let $-\L_{\mathrm{sym}}$ be the nonnegative self-adjoint operator on $L^2_0(\mu)$ associated with $(\Ee_{\mathrm{sym}},\D_0(\Ee_{\mathrm{sym}}))$. Thus $\L_{\mathrm{sym}}$ represents the symmetric part
of $\L$ in the sense of quadratic forms. Formally,
\[
\L_{\mathrm{sym}}=\frac12(\L+\L^*),
\]
where $\L^*$ denotes the adjoint of $\L$ in $L^2(\mu)$. Then
\[
\V=\D\bigl( (-\L_{\mathrm{sym}})^{1/2} \bigr), \qquad
\|u\|_\V=\bigl\| (-\L_{\mathrm{sym}})^{1/2}u \bigr\|_2,
\]
and
\[
\alpha_\V(f)=\|jf\|_{\V^*}
=\|(-\L_{\mathrm{sym}})^{-1/2}f\|_2,
\qquad f\in L^2_0(\mu).
\]
Since $j(L^2_0(\mu))$ is norm dense in $\V^*$, we may identify $\V^*$ with the completion of $L^2_0(\mu)$ under the norm $\|(-\L_{\mathrm{sym}})^{-1/2}\cdot\|_2$.
It remains to control $\beta_\V(f)$. For a bounded measurable function $f :E \to \R$ with $\mu(f) = 0$, \eqref{eq:P} implies
\[
\beta_\V(f) \le  C_P \| f \|_\infty.
\]
Then we obtain the following corollary.

\begin{cor} \label{cor:SP}
Under the Poincar\'e inequality \eqref{eq:P}, let $f: E \to \R$ be a bounded measurable function such that $\mu(f) = 0$. Then for any $t>0$ and $\delta > 0$,
\begin{equation} \label{eq:tail IP'''}
\mathbb P^\mu \left( \left|\frac1t \int_0^{t} f(X_r) \, \d r  \right| \ge \delta  \right) \le 4 \exp \left( - \frac{ t\delta^2}{ 32 \bigl\|(-\L_{\mathrm{sym}})^{-1/2}f \bigr\|_2^2 + 8C_P \| f \|_\infty \delta  } \right).
\end{equation}
The corresponding tail bound for an initial distribution $\nu\ll\mu$ with $\frac{\d\nu}{\d\mu}\in L^p(\mu)$, $1<p\le\infty$, follows from \eqref{eq:tail IP'''} by H\"older's inequality.
\end{cor}

Under an additional Sobolev-Poincar\'e inequality, we can replace the $L^\infty$ bound on $\beta_\V(f)$ by an $L^r$ bound, allowing for unbounded observables. Suppose that, for some $r \in [1, \infty)$ and $C_S > 0$, the following Sobolev-Poincar\'e inequality holds:
\begin{equation}\label{eq:SP}
\|u\|_{2r'} \le C_S \Ee_{\mathrm{sym}}(u,u)^{1/2}, \qquad u\in \D_0(\Ee_{\mathrm{sym}}),
\end{equation}
where $r'$ is the H\"older conjugate of $r$. Since $\mu$ is a probability measure, \eqref{eq:SP} implies \eqref{eq:P} with $C_P=C_S^2$.

Under \eqref{eq:SP}, we have $\V\hookrightarrow L^{2r'}(\mu)$, so Remark~\ref{rem:sobolev} extends the Bernstein bounds to
mean-zero $f\in L^r(\mu)$. For such $f$, we interpret $(-\L_{\mathrm{sym}})^{-1/2}f$ by continuous extension, so that
\[
\|(-\L_{\mathrm{sym}})^{-1/2}f\|_2=\alpha_\V(f).
\]
Together with the bound $\beta_\V(f)\le C_S^2\|f\|_r$, this yields the following corollary.
\begin{cor}
Under the Sobolev-Poincar\'e inequality \eqref{eq:SP}, let $f \in L^r(\mu)$ such that $\mu(f) = 0$. Then for any $t>0$ and $\delta > 0$,
\begin{equation} \label{eq:tail IP''}
\mathbb P^\mu \left( \left|\frac1t \int_0^{t} f(X_s) \, \d s  \right| \ge \delta  \right) \le 4 \exp \left( - \frac{ t\delta^2}{ 32 \bigl\|(-\L_{\mathrm{sym}})^{-1/2}f \bigr\|_2^2 + 8C_S^2 \| f \|_r \delta  } \right).
\end{equation}
The corresponding tail bound for an initial distribution $\nu\ll\mu$ with $\frac{\d\nu}{\d\mu}\in L^p(\mu)$, $1<p\le\infty$, follows from \eqref{eq:tail IP''} by H\"older's inequality.
\end{cor}

As in \cite[Section 3]{GGW}, the estimate reduces to controlling $\mu(fu^2)$ by $\Ee_{\mathrm{sym}}(u,u)$. Here the Sobolev--Poincar\'e inequality gives
$\beta_\V(f)\le C_S^2\|f\|_r$. For diffusion generators, the chain rule yields a more refined bound through the solution of the Poisson equation, as in \cite[Corollary 3.2]{WWZ}. Both estimates thus arise from different bounds on the same quantity $\beta_\V(f)$. In the reversible case, the coefficient $2$ of the asymptotic variance $\sigma^2(f) = 2\alpha_\V(f)^2$ in the Bernstein bounds of \cite{GGW,WWZ} is optimal at the Gaussian scale, as follows from the central limit theorem. While our bound has a larger coefficient, our focus here is on demonstrating the generality of our results rather than on optimizing numerical constants.

We next apply our results to diffusions on compact Riemannian manifolds. Let $(E, g)$ be a $d$-dimensional connected compact Riemannian manifold, possibly with a $C^2$ boundary, and let $V\in C^2(E)$ be normalized so that $\mu(\d x):=\e^{V(x)}\,\d x$ is a probability measure. Let $Z$ be a $C^1$ vector field satisfying
$\mathrm{div}_\mu Z=0$ and tangent to the boundary when $\partial E\ne\emptyset$.
Let $(X_t)_{t\ge0}$ be the diffusion generated by
\[
\L=\Delta+\nabla V\cdot\nabla+Z,
\]
with reflecting boundary when $\partial E\ne\emptyset$. Then $\mu$ is the invariant measure, and the symmetrized Dirichlet form is
\[
\Ee_{\mathrm{sym}}(u,v)
=\int_E\langle\nabla u,\nabla v\rangle_g \,\d\mu,
\qquad u,v\in H^1(E,\mu).
\]
Here $H^1(E,\mu)$ denotes the Sobolev space of functions $u\in L^2(\mu)$ with weak gradient $\nabla u\in L^2(\mu;TE)$. The associated nonnegative self-adjoint operator on $L^2_0(\mu)$ is $-\big(\Delta+\nabla V\cdot\nabla\big)$ with Neumann boundary conditions when $\partial E\ne\emptyset$. Consequently,
\[
\alpha_\V(f)
=\big\|\bigl(-\Delta-\nabla V\cdot\nabla\bigr)^{-1/2}f\big\|_2,
\qquad f\in L^2_0(\mu).
\]
Since $E$ is compact and connected, the Poincar\'e inequality \eqref{eq:P} holds. Moreover, by the Sobolev and Trudinger-Moser inequalities
(see \cite[Chapters 1 and 3]{SC}), the Sobolev-Poincar\'e inequality \eqref{eq:SP} holds with $r=1$ if $d=1$, and with $r=d/2$ if $d\ge3$. In the critical case $d=2$, the Trudinger-Moser inequality and the Orlicz-H\"older inequality yield
\[
\alpha_\V(f) + \beta_\V(f)\le C \|f\|_{L\log L(\mu)}.
\]
Here $\|\cdot\|_{L\log L(\mu)}$ denotes the Luxemburg norm associated with $\Phi(t)=(1+t)\log(1+t)-t$. Since the mean-zero truncations also converge in $L\log L(\mu)$, the preceding approximation argument yields  Bernstein-type inequalities for the additive functionals of $(X_t)_{t\ge0}$ with mean-zero observable
$f$ such that
\begin{equation*}
f \in
\begin{cases}
L^1(\mu) , & \textrm{if } d = 1;\\
L \log L (\mu) , & \textrm{if } d = 2;\\
L^{d/2}(\mu),& \textrm{if } d \ge 3.
\end{cases}
\end{equation*}
See also \cite[Example 3.1, p.~369]{GGW} for the reversible case.

\begin{rem}
The estimate in dimension two is a particular case of an Orlicz embedding argument. Let $\Psi$ be a Young function with convex conjugate $\Psi^*$. Suppose that the following Orlicz-Sobolev-Poincar\'e inequality (see e.g. \cite[Chapter 7]{BGL14}) holds:
\[
\|u^2\|_{L^\Psi(\mu)}
\le C_\Psi \Ee_{\mathrm{sym}}(u, u),
\qquad u\in\V,
\]
then the Orlicz--H\"older inequality (see e.g. \cite{RR}) gives
\[
\beta_\V(f)
\le 2 C_\Psi\|f\|_{L^{\Psi^*}(\mu)},
\]
where both Orlicz norms are Luxemburg norms.
\end{rem}

\section{A corollary of Theorem~\ref{main} and a non-mixing diffusion} \label{sec:non-mixing}

The deviation inequality of Wu \cite{Wu} controls additive functionals through the symmetrized Dirichlet form $\Ee_{\mathrm{sym}}$, which may fail to capture the contribution of the antisymmetric part of the generator to concentration. Motivated by the choices of $\V$ in Section~\ref{section:fin}, we first give a criterion for verifying Assumption~\ref{AS0} in terms of a positive self-adjoint operator and derive the corresponding corollary of Theorem~\ref{main}. We then apply this result to a non-mixing diffusion on the torus. We begin with the following assumption, which leads to Assumption~\ref{AS0}.

\begin{asm}[Resolvent bound associated with a positive self-adjoint operator] \label{AS1}
Let $A$ be a densely defined self-adjoint operator on $L_0^2(\mu)$ such that $A \ge c I$ for some $c>0$, that is,
\[\langle Au,u\rangle_2 \ge c\|u\|_2^2, \qquad u\in\D(A).\]
Suppose that there exists $K > 0$ such that, for every $s>0$ and $g \in L_0^2(\mu)$, we have $R_s g \in \mathscr D (A^{1/2})$ and
\begin{equation} \label{resolvent bound}
\| A^{1/2} R_s g \|_2 \le K \| A^{-1/2} g\|_2.
\end{equation}
\end{asm}

\begin{prp} \label{prop:AS1}
Under Assumption~\ref{AS1}, set
\[
\V = \D(A^{1/2}), \qquad \| u \|_\V = \| A^{1/2} u  \|_2.
\]
Then $\V$ is a Hilbert space, continuously and densely embedded in $L_0^2(\mu)$. Through the canonical embedding $j$, its dual $\V^*$ is the completion of $L_0^2(\mu)$ under the norm $\| A^{-1/2} \cdot \|_2$. Moreover, Assumption~\ref{AS0} holds with the same constant $K$: each $R_s$ admits a unique bounded extension $\widetilde R_s : \V^* \to \V$ satisfying
\[
i \widetilde R_s j = R_s, \qquad \sup_{s>0} \| \widetilde R_s \|_{\V^* \to \V} \le K,
\]
where $i: V \to L_0^2(\mu)$ is the inclusion and $\< jf, g \>_{\V^*, \V} = \< f,ig \>_{2}$.
\end{prp}

\begin{proof}
By spectral calculus, $\V=A^{-1/2}(L^2_0(\mu))$, and $A^{1/2}:\V\to L^2_0(\mu)$ is an isometric isomorphism. Moreover, $\|u\|_2\le c^{-1/2}\|u\|_\V$ for every $u\in\V$, and $\V$ is dense in $L^2_0(\mu)$ since it contains $\D(A)$. Thus $\V$ is a Hilbert space continuously and densely embedded in $L^2_0(\mu)$. For $g \in L_0^2(\mu)$,
\begin{equation} \label{eq:norm}
\| j g \|_{\V^*} = \sup_{\| v \|_2 \le 1} | \langle g, A^{-1/2} v \rangle_2 | = \| A^{-1/2} g \|_2.
\end{equation}
Under the dual identification induced by $A^{1/2}:\V\to L^2_0(\mu)$, the image $j(L^2_0(\mu))$ corresponds to $\V$, which is dense in $L^2_0(\mu)$. Thus $j (L^2_0(\mu))$ is norm dense in $\V^*$. By \eqref{resolvent bound} and \eqref{eq:norm}, the map $j g \mapsto R_s g$ extends uniquely to $\widetilde R_s : \V^* \to \V$ with norm at most $K$, and $i\widetilde R_sj=R_s$ by construction.
\end{proof}

For a bounded measurable $f : E \to \R$ with $\mu(f) = 0$, set
\[
\beta_A (f) = \sup_{u\in\D(A^{1/2}), \, \, \|A^{1/2}u\|_2\le1}|\mu(fu^2)|.
\]
In particular, $\beta_A (f) \le c^{-1} \| f \|_\infty$. Combining Theorem~\ref{main} with Proposition~\ref{prop:AS1}, we obtain the following corollary.

\begin{cor} \label{main'}
Under Assumption \ref{AS1}, let $f: E \to \R$ be a bounded measurable function such that $\mu(f) = 0$. Then for any $t>0$ and $\theta \in \R$ satisfying $4K\beta_A (f)|\theta|\le1$,
\begin{equation*}
\E^\mu \left[ \e^{\theta \int_0^{t} f(X_r) \, \d r}  \right] \le 2 \e^{8K \| A^{-1/2} f \|_2^2 \theta^2 t}.
\end{equation*}
Consequently, for any $\delta > 0$,
\begin{equation*}
\mathbb P^\mu \left( \left|\frac1t \int_0^{t} f(X_r) \, \d r  \right| \ge \delta  \right) \le 4 \exp \left( - \frac{t\delta^2}{8K \big[ 4 \| A^{-1/2} f \|_2^2 + \beta_A (f) \delta \big] } \right).
\end{equation*}
More generally, if $\nu \ll \mu$ and $\frac{\d \nu}{\d \mu} \in L^p(\mu)$ for some $1 < p \le \infty$, then, with $q$ the H\"older conjugate of $p$,
\begin{equation*}
\mathbb P^\nu \left( \left|\frac1t \int_0^{t} f(X_r) \, \d r  \right| \ge \delta  \right) \le 4^{\frac{1}{q}} \left\| \frac{\d \nu}{\d \mu} \right\|_p \exp \left( - \frac{t\delta^2}{8Kq \bigl[ 4 \| A^{-1/2} f \|_2^2 + \beta_A (f) \delta \bigr]} \right).
\end{equation*}
\end{cor}

Corollaries~\ref{cor:BP} and \ref{cor:SP} are recovered from Corollary~\ref{main'} by taking $A=I$ and $A=-\L_{\mathrm{sym}}$, respectively. When the generator $\L$ is normal, i.e. $\L \L^* = \L^* \L$, the modulus of $\L$ provides a natural choice of $A$, which leads to a refined Bernstein-type inequality under an iterated Poincar\'e inequality.

\begin{prp} \label{prop:normal}
Suppose that $\L$ is normal on $L^2(\mu)$ and satisfies the iterated Poincar\'e inequality
\begin{equation} \label{eq:IP'''}
\| \L u \|_2 \ge \eta \| u \|_2, \qquad u \in \mathscr D_0(\L),
\end{equation}
for some $\eta > 0$. With a slight abuse of notation, we denote by $|\L|$ the restriction of $(\L^*\L)^{1/2}$ to $L^2_0(\mu)$. Then $|\L|$ satisfies $|\L| \ge \eta I$. Moreover, for every $g\in L^2_0(\mu)$ and $s>0$, we have $R_sg\in\D(|\L|^{1/2})$ and
\begin{equation} \label{eq:fc}
\bigl\|  |\L|^{1/2} R_s g \bigr\|_2 \le \bigl\| |\L|^{-1/2} g \bigr\|_2, \qquad g \in L_0^2(\mu), \, s > 0.
\end{equation}
\end{prp}

\begin{proof}
The modulus $|\L|$ is nonnegative and self-adjoint, with $\D(|\L|) = \D_0(\L)$ and $\bigl\| |\L| u \bigr\|_2 = \| \L u \|_2$. Thus \eqref{eq:IP'''} and spectral calculus give $|\L| \ge \eta I$.

To apply the spectral theorem for normal operators, we work in the complexification of $L_0^2(\mu)$ and use $\L$ to denote its restriction to this
space. This restriction is still normal, since $L^2_0(\mu)$ is a reducing subspace for $\L$. Let $\sigma(\L) \subset \C$ be the spectrum of $\L$ and $\{E_\lambda : \lambda \in \sigma(\L) \}$ the spectral family corresponding to $\L$. Since $\L$ generates a contraction semigroup, $\sigma(\L) \subset \{\lambda : \mathrm{Re} \, \lambda \le 0  \}$. Therefore, for every $\lambda \in \sigma(\L)$ and $s>0$,
\[
|s - \lambda|^2 = |\lambda|^2 + s^2 - 2s \mathrm{Re} \, \lambda \ge |\lambda|^2.
\]
On the other hand, \eqref{eq:IP'''} implies that $\sigma(\L) \subset \{\lambda : |\lambda| \ge \eta  \}$. Then it follows from functional calculus that, for any $g \in L_0^2(\mu)$ and $s > 0$,
\[
\bigl\|  |\L|^{1/2} R_s g \bigr\|_2^2 = \int_{\sigma(\L)} \frac{|\lambda|}{|s - \lambda|^2} \, \d \| E_\lambda g\|_2^2 \le \int_{\sigma(\L)} \frac{1}{|\lambda|} \, \d \| E_\lambda g\|_2^2 = \bigl\| |\L|^{-1/2} g \bigr\|_2^2 < \infty.
\]
This proves that $R_s g \in \D(|\L|^{1/2})$ and gives \eqref{eq:fc}. The proof is complete.
\end{proof}

Under the assumptions of Proposition~\ref{prop:normal}, we can take $A = |\L|$ and apply Corollary~\ref{main'} with $K = 1$. Given any bounded mean-zero $f$, for $t>0$ and $\delta>0$,
\begin{equation} \label{eq:tail'}
\mathbb P^\mu \left( \left|\frac1t \int_0^{t} f(X_r) \, \d r  \right| \ge \delta  \right) \le 4 \exp \left( - \frac{t\delta^2}{32 \bigl\| |\L|^{-1/2} f \bigr\|_2^2 + 8 \beta_{|\L|} (f) \delta  } \right),
\end{equation}
where
\[
\beta_{|\L|} (f) = \sup_{u\in\D(|\L|^{1/2}), \, \, \| |\L|^{1/2}u\|_2\le1}|\mu(fu^2)| \le \eta^{-1} \| f \|_\infty.
\]
Unlike the symmetrized Dirichlet form, the modulus $|\L|$
also reflects the antisymmetric part of the generator.
The following example shows how this yields Bernstein-type
inequalities for a non-mixing diffusion.

Let $\mathbb T = \R/(2\pi \Z)$ be the torus, and let $\mu$ be the normalized Lebesgue measure on $\mathbb T^2$. Consider the process on $\mathbb T^2$
\[
X_t = X_0 + (t, \sqrt{2} W_t)  \pmod {2\pi},
\]
where $(W_t)_{t\ge0}$ is a standard Brownian motion independent of $X_0$. Its generator, in coordinates $(x,y)$, is
\[
\L=\partial_x+\partial_y^2.
\]
Integration by parts shows that the associated closed symmetrized Dirichlet form is given by
\[
\Ee_{\mathrm{sym}}(u, u) = \| \partial_y u \|_2^2, \qquad \mathscr D(\Ee_{\mathrm{sym}}) = L^2(\mathbb T_x ; H^1(\mathbb T_y)).
\]
The measure $\mu$ is invariant. Let ${\mathrm H}_t= \e^{t\partial_y^2}$ denote the heat semigroup acting on the second variable. Then, for every $u\in L^1(\mu)$
and $t\ge0$,
\[
P_tu(x,y) = {\mathrm H}_t\bigl[u(x+t,\cdot)\bigr](y), \qquad \mu\text{-a.e. }(x,y)\in\mathbb T^2.
\]
To see that the process is not mixing, consider $u(x,y)=f(x)$ with $f\in L^1(\mathbb T)$. Then $P_tu(x,y)=f(x+t)$. In particular, for $u(x,y)=\cos x$, we have
\[
\mu(uP_tu)=\frac12\cos t,
\]
which does not converge to $\mu(u)^2=0$ as $t\to\infty$. Moreover, $\Ee_{\mathrm{sym}}$ vanishes on functions depending only on $x$, so Wu's inequality yields only a trivial bound for $u=\cos x$ when $0<\delta<1$.

The adjoint of $\L$ is
\[
\L^*=-\partial_x+\partial_y^2,
\]
with the same domain as $\L$, and $\L$ is normal. On the complexification of $L^2_0(\mu)$, the functions
\[
e_{k,\ell}(x,y):=\e^{(kx+\ell y)\i}, \qquad (k,\ell)\in\mathbb Z^2\setminus\{(0,0)\},
\]
form an orthonormal basis, and
\[
\L e_{k,\ell}=(k \i -\ell^2)e_{k,\ell},
\qquad
\L^*e_{k,\ell}=(-k \i -\ell^2)e_{k,\ell}.
\]
Thus $\L$ is normal, and $|\L|$ satisfies
\[
|\L|e_{k,\ell}=\rho_{k,\ell}e_{k,\ell}, \qquad \rho_{k,\ell}:=\sqrt{k^2+\ell^4}.
\]
Writing $\widehat u(k,\ell)$ for the Fourier coefficients of $u$,
we have
\[
\D(|\L|)
=
\bigg\{
u\in L^2_0(\mu):
\sum_{(k,\ell)\ne(0,0)}
(k^2+\ell^4)|\widehat u(k,\ell)|^2<\infty
\bigg\}.
\]
Since $\rho_{k,\ell}\ge1$ for every $(k,\ell)\ne(0,0)$,
we have $|\L|\ge I$. Hence
\[
\|\L u\|_2=\bigl\||\L|u\bigr\|_2\ge\|u\|_2,
\qquad u\in\D_0(\L),
\]
so Proposition~\ref{prop:normal} applies with $\eta=1$. By the spectral representation of $|\L|$,
\[
\D(|\L|^{1/2})
= \bigg\{
u\in L^2_0(\mu): \sum_{(k,\ell)\ne(0,0)} (k^2+\ell^4)^{1/2}|\widehat u(k,\ell)|^2<\infty
\bigg\},
\]
and
\begin{equation} \label{eq:norm'}
\big\||\L|^{1/2}u\big\|_2^2
=
\sum_{(k,\ell)\ne(0,0)}
(k^2+\ell^4)^{1/2}|\widehat u(k,\ell)|^2.
\end{equation}
We next show that $\D(|\L|^{1/2}) \hookrightarrow L^6(\mu)$. Let $J_x = (I -\partial_x^2)^{1/2}$ and $J_y = (I -\partial_y^2)^{1/2}$. The one-dimensional Sobolev embedding $H^{1/3}(\mathbb T) \hookrightarrow L^6(\mathbb T)$, together with its Hilbert-valued version, gives
\begin{equation} \label{eq:embed 6}
\| u \|_6 = \| u \|_{L_x^6 L_y^6} \le C \| J_y^{1/3} u\|_{L_x^6 L_y^2} \le C' \| J_x^{1/3} J_y^{1/3} u\|_{L_x^2 L_y^2},
\end{equation}
for every smooth mean-zero function $u$ on $\mathbb T^2$. On the other hand, by Parseval's identity,
\begin{equation} \label{eq:norm''}
\| J_x^{1/3} J_y^{1/3} u\|_{L_x^2 L_y^2}^2 = \sum_{k, \ell \in \Z} (1+k^2)^{1/3} (1+\ell^2)^{1/3} |\widehat u(k,\ell)|^2.
\end{equation}
By Young's inequality, for $(k,\ell)\ne(0,0)$,
\begin{align*}
(1+k^2)^{1/3} (1+\ell^2)^{1/3} &\le (1 + |k|)^{2/3} (1+\ell^2)^{1/3} \le \frac23(1 + |k|) + \frac13 (1+\ell^2)\\
& \le 2(k^2+\ell^4)^{1/2},
\end{align*}
Combining \eqref{eq:norm'} and \eqref{eq:norm''}, we obtain
\[
\| J_x^{1/3} J_y^{1/3} u\|_{L_x^2 L_y^2}^2 \le 2 \big\||\L|^{1/2}u\big\|_2^2.
\]
Together with \eqref{eq:embed 6}, this yields
\[
\|u\|_6 \le 2C' \big\||\L|^{1/2}u\big\|_2.
\]
By density of mean-zero smooth functions in $\D(|\L|^{1/2})$, this bound extends to every $u\in\D(|\L|^{1/2})$. Consequently, combining \eqref{eq:tail'} with Remark~\ref{rem:sobolev} and the preceding embedding, we obtain the following Bernstein-type inequality for $(X_t)_{t \ge 0}$: there exists $C > 0$ such that, for every mean-zero $f \in L^{3/2}(\mu)$, $t>0$ and $\delta>0$,
\begin{equation} \label{eq:tail''}
\mathbb P^\mu \left( \left|\frac1t \int_0^{t} f(X_r) \, \d r  \right| \ge \delta  \right) \le 4 \exp \left( - \frac{t\delta^2}{32 \bigl\| |\L|^{-1/2} f \bigr\|_2^2 + C \| f \|_{3/2} \delta  } \right),
\end{equation}
Here $|\L|^{-1/2}f$ is understood by continuous extension from $L^2_0(\mu)$ to the mean-zero subspace of $L^{3/2}(\mu)$.

\paragraph{Data Availability Statements.} Data sharing not applicable to this article as no datasets were generated or analysed during the current study.

\paragraph{Acknowledgement.} The author used ChatGPT-6 Astra during the development of this work for exploratory discussions of proof strategies and for language assistance. All mathematical arguments were independently verified and formalized by the author.

\end{document}